\documentclass{amsart}
\usepackage[utf8]{inputenc}
\usepackage{lmodern}
\usepackage[margin=3cm]{geometry}
\usepackage{extdash} 
\usepackage{amsmath,amsthm,amsfonts,amssymb, microtype}
\usepackage{xr-hyper}
\usepackage[colorlinks=true]{hyperref}
\usepackage{verbatim}
\usepackage{pdfsync}
\usepackage{stmaryrd}
\usepackage{marvosym}
\usepackage{colonequals}
\usepackage{mathrsfs}
\usepackage{bbold}

\usepackage{mathtools}
\usepackage{enumitem}
\usepackage{tikz}
\usetikzlibrary{arrows,calc,positioning}

\usetikzlibrary{backgrounds}

\usepackage{extdash}
\usepackage{amsthm,amsfonts,amssymb,amsmath, microtype}
\usepackage{xr-hyper}
\usepackage[colorlinks=true]{hyperref}
\usepackage{verbatim}
\usepackage{pdfsync}
\usepackage{comment} 

\usetikzlibrary{circuits.logic.IEC}
\usetikzlibrary{cd}
\usepackage{xcolor}
\usetikzlibrary{shapes}
\usetikzlibrary{decorations.markings}

\tikzstyle heightone=[scale=.7,shift={(0,-.3)}]
\tikzstyle heightones=[scale=.8,xscale=.35,shift={(0,.1)}]
\tikzstyle heightoneonehalf=[scale=.9,shift={(0,-.2)}]
\tikzstyle heighttwo=[scale=.9,shift={(0,-.4)}]
\tikzstyle heighttwos=[scale=.5,xscale=.6,shift={(0,-.1)}]
\tikzstyle heightthree=[scale=.6,shift={(0,-.9)}]
\tikzstyle heightthrees=[scale=.4,xscale=.7,shift={(0,-.2)}]

\tikzstyle arrowstyle=[blue,semitransparent,scale=2]

\tikzstyle basiclabel=[draw=none,fill=none,shape=rectangle,inner sep=2pt,scale=.8]
\tikzstyle leftlabel=[basiclabel,anchor=east]
\tikzstyle rightlabel=[basiclabel,anchor=west]
\tikzstyle bottomlabel=[basiclabel,anchor=north]
\tikzstyle toplabel=[basiclabel,anchor=south]

\tikzstyle vertex=[circle,draw,fill=black,inner sep=1pt]
\tikzstyle ciliation=[circle,draw=none,fill=red,inner sep=1pt,semitransparent]
\tikzstyle ciliatednode=[vertex,pin={[pin distance=1mm,pin edge={semitransparent,red},ciliation]#1:{}}]

\tikzstyle vector=[black,thick,rectangle,draw=gray!50!yellow,top color=yellow!30,bottom color=black!10,scale=.8,inner sep=2pt]
\tikzstyle small vector=[vector,scale=.8]
\tikzstyle plain vector=[rectangle,draw=none,fill=white,scale=.7]

\tikzstyle my signal=[black,thick,signal,signal pointer angle=120,draw=blue!50,top color=blue!20,bottom color=black!10,scale=.8,inner sep=2pt]
\tikzstyle matrix=[my signal,signal from=south,signal to=north]
\tikzstyle reverse matrix=[my signal,signal from=north,signal to=south]

\tikzstyle small matrix=[matrix,scale=.7]
\tikzstyle reverse small matrix=[reverse matrix,scale=.7]
\tikzstyle matrix on edge=[small matrix,sloped,rotate=-90]
\tikzstyle reverse matrix on edge=[small matrix,sloped,rotate=90]

\tikzstyle trivalent=[very thick]
\tikzstyle dotdotdot=[decorate,decoration={markings,
    mark=at position .3 with{\node{.};},
    mark=at position .5 with {\node{.};},
    mark=at position .7 with {\node{.};}}]

\tikzstyle wavyup=[out=90,in=-90]
\tikzstyle wavydown=[out=-90,in=90]

\tikzstyle symmetrizer=[rectangle,fill=gray!10,draw=black]
\tikzstyle permutation=[symmetrizer]
\tikzstyle antisymmetrizer=[rectangle,fill=black,draw=black]
\tikzstyle symlabel=[draw=none,fill=none,black,scale=.8]
\tikzstyle asymlabel=[draw=none,fill=none,white,scale=.8]

\newcommand{\Z}{\mathbb{Z}}

\newcommand{\N}{\mathbb{N}}
\newcommand{\Sphere}{\mathbb{S}}
\newcommand{\DDD}{\mathbb{D}}

\newcommand{\K}{\mathrm{K}}

\newcommand{\id}{\mathrm{Id}}

\newcommand{\modules}{\mathrm{-mod}}

\newcommand{\A}{\mathcal{A}}
\newcommand{\C}{\mathcal{C}}

\newcommand{\1}{\mathbb{1}}

\newcommand{\Top}{\underline{\mathrm{Top}}}

\newcommand{\D}{\mathcal{D}}

\newcommand{\kk}{\Bbbk}

\newcommand{\nocontentsline}[3]{}
\newcommand{\tocless}[2]{\bgroup\let\addcontentsline=\nocontentsline#1{#2}\egroup}

\newtheorem{theorem}{Theorem}[section]
\newtheorem{proposition}[theorem]{Proposition}
\newtheorem{lemma}[theorem]{Lemma}

\newtheorem*{theorem*}{Theorem}

\theoremstyle{definition}
\newtheorem{definition}[theorem]{Definition}
\newtheorem{lemma/definition}[theorem]{Definition/Lemma}
\newtheorem{example}[theorem]{Example}
\newtheorem{observation}[theorem]{Observation}

\newtheorem{remark}[theorem]{Remark}

\begin{document}

\title{Algebraic $K_0$ for unpointed Categories}
\author{Felix K\"ung}
\begin{abstract}
We construct a natural generalization of the Grothendieck group $\mathrm{K}_0$ to the case of possibly unpointed categories admitting pushouts by using the concept of heaps recently
introduced by Brezinzki. In case of a monoidal category, the defined K0 is
shown to be a truss. It is shown that the construction generalizes the classical $\K_0$ of an abelian category as the group retract along the isomorphism class of the zero object. We finish by applying this construction to construct the integers with addition and multiplication as the decategorification of finite sets and show that in this $\K_0\left(\underline{\mathrm{Top}}\right)$ one can identify a CW-complex with the iterated product of its cells.

\end{abstract} 
\maketitle
\let\thefootnote\relax\footnotetext{The author is a Postdoctoral Research Fellow at the Universit\'e Libre de Bruxelles. The research in this work has been supported by the MIS (MIS/BEJ - F.4545.21) Grant from the FNRS and an ACR Grant from the Universit\'e Libre de Bruxelles.}
\section*{Introduction}

The Grothendieck group $\K_0\left(\A\right)$ of an abelian category $\A$ is one of the most studied and used invariants of an abelian category \cite{Weibel1994,Milnor1971,Brown1973,Quillen1975,Rosenberg1994} and has been generalized to the setting of derived and triangulated categories \cite{Neeman1997} and equivariant objects \cite{Schwede2022,Weibel2013}. In the setting of homotopy categories the general idea is to replace the short exact sequences in the definition of $\K_0\left(\A\right)$ with triangles \cite{Neeman1997,Neeman1998,Neeman1998a,Neeman1999}, respectively an exact structure \cite{Waldhausen1985}, both of which can often be thought about as shadows of short exact sequences. Another approach of generalization replaces the abelian structure with an additive structure to generate the split Grothendieck Group $\K_0^\oplus\left(\A\right)$ \cite{Swan1968}. However, this approach remembers significantly less structure and only recovers the classical definition if the category $\A$ is semisimple abelian. 

A useful feature of these groups is that they turn into rings if the abelian category admits the structure of a monoidal category which respects the abelian structure. These have been used a lot over the years in order to categorify algebraic structures like Hecke algebras \cite{Soergel2007}, representations \cite{Stroppel2005, Libedinsky2022}, quantum groups along roots of unity \cite{Elias2016} and many more.

However, the assumption on a category to be abelian respectively pointed is strong, so it would be useful to generalize the construction of the Grothendieck group to a more general setting. In this paper we present a generalization of $\K_0$ to the case of a category admitting pushouts. One natural generalization of a group is a heap \cite{Brzezinski2020,Brzezinski2022a}, which is also compatible with the natural generalization of a ring, a truss \cite{Brzezinski2019,Brzezinski2022}. In this paper we take the point of view of so-called heaps to construct an analog of the Grothendieck group for unpointed categories. This in particular allows us to consider $\K_0$ of the category of topological spaces, which does not admit a zero object, as the initial object $\{*\}$ and final object $\{\}$ do not coincide.

All in all we show in this paper that $\K_0\left(\A\right)$ admits a natural generalization to non-abelian categories by passing from groups to heaps, and rings to trusses. In particular the structure of it being group, respectively a ring can be seen as a shadow of having a zero object, respectively being pointed. 

In Definition~\ref{Definition heapy K0} we define the Grothendieck heap of a category admitting pushout (respective pullback) we then prove in Theorem~\ref{Lemma Heapy K0 recovers classical K0} that this recovers the classical definition of $\K_0\left(\A\right)$ for an abelian category and in Theorem~\ref{Theorem Truss structure on K0} that it admits naturally the structure of a truss if the starting category was monoidal. Again this construction recovers the original definition of the ring structure on $\K_0\left(\A\right)$ for a monoidal abelian category $\A$. 

Finally we illustrate the properties of our construction on a few examples: First in Example~\ref{Example Sets} we prove that the integers can be recovered as $\K_0\left(\underline{\mathrm{Set}}\right)$, where $\underline{\mathrm{Set}}$ denotes the category of finite sets. Then we finish with Example~\ref{Example Top} by observing that for topological spaces we have in $\K_0\left(\underline{\mathrm{Top}}\right)$ that CW-complexes can be identified with the iterated product of their cells.

\section*{Acknowledgements}
I would like to thank P. Sarraco, J. Vercruysse and W. Hautekiet for enlightening comments and introducing me to the subject of heaps respectively trusses. I am deeply grateful for G. Janssens and \v{S}. \v{S}penko, for the awesome comments, hints and support while writing this paper. Finally I would like to show my gratitude to L. Gilbert whose talk on her bachelor thesis made me realise that heaps behave like decategorified pushouts.

I also want to thank the referee for the helpful comments and suggestions.

\section{Preliminaries}

In this section we collect some basic notions about heaps and the classical algebraic $\K$-groups $\K_0$ and $\K_0^\oplus$.

\subsection{Heaps and their morphisms}

We start by recalling a few notions about heaps and trusses, for a more in depth treaties we refer to \cite{Brzezinski2020}.

\begin{definition}
A heap is a pair $ \left(H,[\_,\_,\_] \right) $ where $[\_,\_,\_]:H\times H \times H\to H$ is such that:
\begin{align*}
[a,b,[c,d,e]]&=[[a,b,c],d,e]&\text{(Associativity)}\\
[x,x,y]= &y=[y,x,x]&\text{(Identity)}.
\end{align*} 
\end{definition}

In particular we may write without loss of generality $[x_1,...,x_n]$ for iterated ternary products as long as $n$ is odd. 
We also set $[x]=x$ for convenience.

\begin{example}\label{example Group}
Let $G$ be a group. Then $G$ defines canonically a heap via the operation  $\left[x,y,z\right]:=x  y^{-1}z$.
\end{example}

\begin{definition} Let $H, H'$ be two heaps, a map $\varphi: H \to H'$ is a morphism of heaps if $$[\varphi\left(x\right),\varphi\left(y\right),\varphi\left(z\right)]=\varphi\left([x,y,z]\right).$$
\end{definition}

Even though the following definition of a free abelian heap might seem complicated it mimics the construction of free abelian groups and inherits the same universal property.

\begin{definition}
Let $X$ be a set. Then the free abelian heap over $X$ is given by the set of equivalence classes of words $$\left\{\left\{x_{\sigma 1},y_{\sigma ' 1},x_{\sigma 2},...y_{\sigma ' n-1}, x_{\sigma n}\right\}_{\sigma\in S^n,\sigma'\in S^{n-1}}\mid x_\sigma{i} \neq y_\sigma{i} \right\}$$
where $\sigma \in S^n$ and $\sigma \in S^{n-1}$. With heap structure given by concatenating the sequence and then reducing.
\end{definition}

Together with the notion of a normal subheap this allows us to effectively define heaps via free abelian heaps and relations.

\begin{definition}
Let $H$ be a heap. A subheap $S\subset H$ is normal if we have for every $x\in H$ and $e,s\in S$ a $t\in S$ such that $\left[\left[x,e,s\right],x,e\right]\in S$.
\end{definition}

\begin{proposition}
Let $H$ be a heap and $S\subset H$ a subheap. Then $H/S$ is a heap with operation $$\left[\overline{x},\overline{y},\overline{z}\right]:=\overline{\left[x,y,z\right]}.$$
\end{proposition}

\begin{definition}\label{Definition truss} A truss $\left(T,\left[\_,\_,\_\right],\cdot \right)$ is a heap $\left(T , \left[\_,\_,\_\right]\right)$ together with a multiplication $\cdot: T\times T\to T$ such that $$\left[w\cdot x,w \cdot y, w \cdot z\right]=w\cdot\left[x,y,z\right]\qquad w,x,y,z \in T.$$
\end{definition}

As usual we will write mostly $xy$ for $x\cdot y$.

We will use the following construction to freely pass between heaps with a natural neutral element and groups.

\begin{definition}
Let $H$ be a heap and $e \in H$. Then the retract of $H$ along $e$ is the group $\left(H,+_e\right)$ defined to have the same underlying set as $H$, but with addition $a+_e b := [a,e,b]$.
\end{definition}

In the group $\left(H,+_e\right)$ we have $e$ as neutral element. Furthermore one can compute for $x \in H$
\begin{align*}
[x,e,[e,x,e]]&=[[x,e,e],x,e] &\text{(Associativity)}\\
&=[x,x,e] &\text{(Identity)}\\
&=e &\text{(Identity)}
\end{align*}
and so $[e,x,e]$ is the inverse of $x$.

\begin{lemma}\label{Proposition passing between groups and heaps with dedicated 0}
Let $H$ be a heap and let $e\in H$. Then we have that $[a,b,c]=a -_e b +_e c $. In particular we can construct the original heap from its retract.
\end{lemma}

\begin{lemma}\label{Lemma induced group morphism}
Let $\varphi: H \to H'$ be a morphism of heaps such that $\varphi\left(e\right)=e'$. Then $\varphi$ induces a canonical group morphism $$\left(H,+_e\right)\to \left(H',+_{e'}\right).$$
\end{lemma}

\begin{definition}
Let $\left(T,\left[\_,\_,\_\right]_T,\cdot_T \right)$ and $\left(T',\left[\_,\_,\_\right]_{T'},\cdot_{T'} \right)$ be two trusses. A morphism of heaps $T\to T'$ is then a morphism of trusses if additionally 
$$\varphi\left(xy\right)=\varphi\left(x\right)\varphi \left(y\right).$$ \end{definition}

\begin{definition}
A sub-heap $S$ of a truss $\left(T,\left[\_,\_,\_\right],\cdot\right)$ is an ideal if we have for all $x\in T$ and $s\in S$
$$xs \in S  \quad sx \in S$$ 
\end{definition}

In general in order to get a canonical truss structure on the quotient of a truss by a subheap one needs the subheap to be a so called paragon. However, as we will verify that a normal subheap is actually an ideal we may consider quotients along ideals.

\begin{lemma}
Let $S\subset T$ be an ideal of a truss, then we have a canonical truss structure on $T/S$.
\end{lemma}

\subsection{Algebraic $\K_0$}

We now recall the following definitions and constructions in order to later imitate them on the non-abelian level. The notions we cover here are very classical and there are countless references for them, for instance \cite{Swan1968,Rosenberg1994}.

\begin{definition}
Let $\A$ be an (essentially small) abelian category, then the Grothendieck group $\K_0\left(\A\right)$ is the free abelian group of isomorphism classes of objects in $\A$ with sums induced by short exact sequences. More precisely:
$$\K_0\left(\A\right):=\langle \overline{A}\mid \overline{B}=\overline{A}+ \overline{C} \text{ if there exists a short exact sequence } 0 \to A \hookrightarrow B \twoheadrightarrow C\to 0 \rangle$$
\end{definition}

\begin{proposition}
Let $F: \A \to \A' $ be an exact functor of abelian categories (i.e. sending short exact sequences to short exact sequences). Then $F$ induces a morphism of groups
\begin{align*}
F:\K_0\left(\A\right) &\to \K_0\left(\A'\right)\\
\overline{A}&\mapsto \overline{F \left(A\right)}.
\end{align*}
\end{proposition}

In absence of exact sequences one still has the following notion for additive categories.
\begin{definition}
Let $\A$ be an (essentially small) additive category, then the split Grothendieck group $\K_0^\oplus\left(\A\right)$ is defined to be the free abelian group on isomorphism classes of $\A$ with relations induced by $\oplus$. More precisely:
$$\K_0^\oplus\left(\A\right):=\langle \overline{A}\mid \overline{A\oplus B}=\overline{A}+\overline{B}\rangle,$$
where we denote by $\overline{A}$ the isomorphism class of $A$.
\end{definition}

\begin{proposition}
Let $\A \to \A'$ be an additive functor of additive categories (i.e. $F$ preserves direct sums). Then $F$ induces a morphism of groups:
\begin{align*}
F:\K_0^\oplus\left(\A\right)&\to \K_0^\oplus\left(\A'\right)\\
\overline{A}&\mapsto \overline{\left(A\right)}
\end{align*}
\end{proposition}


\begin{example} We give examples for the $\K_0$ respectively $\K_0^\oplus$ of some well known categories:
\begin{itemize}
\item Let $\kk$ be a field, then we have that $\K_0\left(\kk\modules\right)\cong \mathbb{Z}\cong \K_0^\oplus \left(\kk\
\modules\right)$, as modules over a field are vector spaces and these are uniquely defined up to isomorphism by their dimensions. Furthermore the dimension of vector spaces is additive with respect to the direct sum and additive with respect to short exact sequences, so the defining relations literally coincide.
\item We have $\K_0\left(\mathbb{Z}\modules \right)\cong \mathbb{Z}$ but a non-identity projection $$\K_0^\oplus\left(\mathbb{Z}\modules\right)\twoheadrightarrow \K_0\left(\mathbb{Z}\modules \right)\cong \mathbb{Z}.$$ 
The reason for this is that the object $\mathbb{Z}\otimes \left(\mathbb{Z}/n\mathbb{Z}\right)$ is not isomorphic to $\mathbb{Z}$ and cannot be built from it by direct sums. However, we have a non-split short exact sequence
$$0\to \mathbb{Z}\xhookrightarrow{n\cdot \_} \mathbb{Z}\twoheadrightarrow \mathbb{Z}/n\mathbb{Z}\to 0.$$
So $\overline{\mathbb{Z}/n\mathbb{Z}}+\overline{\mathbb{Z}}=\overline{\mathbb{Z}}$ and in particular $\overline{\mathbb{Z}/n\mathbb{Z}}=0$.
\end{itemize}
\end{example}

\begin{remark}
In the examples above it is important to consider the categories of finitely generated modules, respectively vector spaces, as otherwise the whole group collapses due to something called the "Eilenberg Swindle".

For example in the case of infinite dimensional vector spaces, i.e. $\K_0\left(\kk\mathrm{-Mod}\right)$, there exists a vector space $A:=\bigoplus_{i\in \N}\kk$ such that for all finite dimensional vector spaces $V$ we have $V\oplus A\cong A$. This means that in $\K_0\left(\kk-\mathrm{Mod}\right)$ we have $\overline{V}+\overline{A}=\overline{A}$ and so $\overline{V}=0$. Now if $V$ is infinite dimensional, then $V\oplus V\cong V$ and so we have analogously $\overline{V}=0$ too and $$\K_0\left(\kk\mathrm{-Mod}\right)\cong 0.$$
\end{remark}

\begin{observation} Direct sums and exact sequences can be considered as special cases of pushouts, respectively pullbacks. In particular we consider for a pointed category admitting coproducts the direct sum as a pushout along the zero object.
\end{observation}

For the convenience of the reader and in order to present the general arguments we include the proofs of the following well known results.

\begin{lemma}\label{Lemma split K0 surjects onto ordinary K0}
Let $\A$ be an abelian category, then we have a natural surjection 
$$\K_0^\oplus\left(\A\right)\twoheadrightarrow \K_0\left(\A\right).$$
\begin{proof}
We will denote the isomorphism class of $X$ in $\K_0^\oplus\left(\A\right)$ by $\overline{X}^\oplus$.
We have a canonical map on generators by $\overline{X}^\oplus\mapsto \overline{X}$. So we need to verify that the relations of $\K_0^\oplus\left(\A\right)$ are preserved. In particular we need to show that in $\K_0\left(\A\right)$ we have the relation $\overline{X\oplus Y}=\overline{X}+\overline{Y}$. To recover this relation consider the split short exact sequence 
$$0\to A \hookrightarrow A\oplus B  \twoheadrightarrow B \to 0$$ which gives rise to the generating relation $\overline{A\oplus B}^\oplus= \overline{A}+\overline{B},$ and so the identity functor induces the claimed projection.
\end{proof}
\end{lemma}


\begin{lemma}
Let $\A$ be an abelian category with a monoidal structure such that $-\otimes -$ is exact. Then $\K_0\left(\A\right)$ has naturally the structure of a ring with product given by $$\overline{A}\overline{B}:=\overline{A\otimes B}.$$
\begin{proof}
As we have a natural multiplication on $\left\{\overline{A}\right\}_{A\in \A}$ we get a natural multiplication on 
$$R:=\langle \overline{A}\rangle_{A\in\A}.$$ 
This makes $\langle \overline{A}\rangle_{A\in\A}$ into an ring with unit given by $\overline{\1}$ for $\1$ the monoidal unit. Now we need to prove that the subgroup 
$$I = \langle\overline{A} +\overline{C}-\overline{B}\mid 0\to A\hookrightarrow B \twoheadrightarrow C \to 0 \text{ short exact sequence}\rangle$$
 is a two-sided ideal in $R$.

This we can check by distributivity on generators. So consider a generator $\left(\overline{A}+\overline{C}-\overline{B}\right)$ induced by the short exact sequence $0\to A \hookrightarrow B\twoheadrightarrow C \to 0$ and $\overline{X}\in \K_0\left(\A\right)$ with $X\in \A$. Then we have by exactness of $-\otimes X$ that 
$$0\to A\otimes X \hookrightarrow B\otimes X\twoheadrightarrow C\otimes X \to 0$$
is exact and so 
\begin{align*}\left(\overline{A\otimes X}+\overline{C\otimes X}-\overline{B\otimes X}\right)&=\left(\overline{A}\;\overline{X}+\overline{C}\;\overline{X}-\overline{B}\;\overline{X}\right)\\
&=\left(\overline{A}+\overline{C}-\overline{B}\right)\overline{X}
\end{align*} 
which means that $I$ is a right ideal. However, as we only used exactness of $-\otimes-$ we can apply the same argument for left multiplication and so $I$ becomes a two sided ideal and $$\K_0\left(\A\right)=R/I$$
has a natural multiplication induced by $\overline{A}\;\overline{B}=\overline {A\otimes B}$. 
\end{proof}
\end{lemma}

\begin{remark}
Often one can generalize the above result to considering right exact monoidal structures. However, this usually is done by representing the class of every object in $\K_0\left(\A\right)$ as sum of flat objects and then applying $-\otimes-$ there. In particular we would need some kind of assumption on having enough flat objects, which we want to avoid at the moment.
\end{remark}

\begin{proposition}
Let $F: \A \to \A'$ be an exact monoidal functor between abelian categories, i.e. $F\left(\1\right)\cong \1$ and $F\left(A\otimes B\right)\cong F\left(A\right)\otimes F\left(B\right)$. Then $F$ induces a morphism of rings 
\begin{align*}
\K_0\left(\A\right)&\to \K_0\left(\A'\right)\\
\K_0\left(\A\right)&\mapsto \K_0\left(\A'\right)
\end{align*}
\begin{proof}
We know that $F$ induces a group morphism $\K_0\left(\A\right)\to \K_0\left(\A'\right)$, so we need to check that it preserves the product. But this follows immediately as $F$ is monoidal and so 
\begin{align*}
&F \overline{\1}=\overline{F\left(\1\right)}=\overline{\1}\\
&F\left(\overline{A}\;\overline{B}\right)=F\left(\overline{A\otimes B}\right)=\overline{F\left(A \otimes B\right)}=\overline{F\left(A\right)\otimes F\left(B\right)}=\overline{F\left(A\right)}\overline{F\left(B\right)}=F\left(\overline{A}\right)F\left(\overline{B}\right)
\end{align*}
which proves that $F$ is a ring morphism.
\end{proof}
\end{proposition}

\section{$\K_0$ for categories admitting pushouts}

We will for ease of reading restrict to the case of categories admitting pushouts. The case for pullbacks is analagous after passing to the opposite category, see Remark~\ref{Remark how to dualize pushouts}.

\subsection{The Grothendieck heap and relation to the classical setting}

We start by defining the Grothendieck heap on a category admitting pushouts and monomorphisms and verify that it indeed is a generalization of the classical $\K_0\left(\A\right)$ for an abelian category. However, first we give the following observation showing that heap relations arise naturally in the isomorphism classes of objects.

\begin{observation}
Pushouts satisfy categorical incarnations of the heap relations, i.e.\ if we place objects in the corners of iterated pushouts, then this construction satisfies (Associativity) and (Identity) up to natural isomorphisms.
\begin{proof}
(Associativity): By pacing of pushouts we have that up to unique isomorphisms iterated pushouts do not depend on the sequence of pushing out. In particular this is a categorical incarnation of (Associativity).

(Identity): Consider the two diagrams
$$\tikz[heighttwo,xscale=2,yscale=2,baseline]{
\node (X1) at (0,1) {$X$};
\node (X2) at (0,0) {$X$};
\node (Y1) at (1,1) {$Y$};
\node (Y2) at (1,0) {$Y$};
\node at (2,0.5){and};
\node (X3) at (3,1) {$X$};
\node (X4) at (3,0) {$Y$};
\node (Y3) at (4,1) {$X$};
\node (Y4) at (4,0) {$Y.$};
\node (po)at (0.5,0.5) {$\ulcorner$};
\node (po)at (3.5,0.5) {$\ulcorner$};

\draw[->]
(X1) edge node[left] {$\id$} (X2)
(X1)edge node [above] {$\alpha$} (Y1)
(X2)edge node [above] {$\alpha$} (Y2)
(Y1)edge node [right] {$\id$} (Y2)
(X3) edge node[left] {$\beta$} (X4)
(X3)edge node [above] {$\id$} (Y3)
(Y3)edge node [right] {$\beta$} (Y4)
(X4)edge node [above] {$\id$} (Y4)
;
}
$$
inducing the (Identity) relation.

So pushouts naturally induce heap relations on the isomorphism classes of objects.
\end{proof}
\end{observation}

\begin{definition}\label{Definition heapy K0}
Let $\C$ be an (essentially small) category admitting pushouts. Then we define the Grothendieck heap of $\C$ to be the free abelian heap generated by isomorphism classes of objects in $\C$ with relations induced by pushouts admitting at least one monomorphic leg. More precisely
$$\K_0\left(\C\right):=\langle \overline{X}\mid [\overline{X},\overline{Y},\overline{Z}]= \overline{X\sqcup_Y Z} \text{ for }Y\hookrightarrow X \text{ or } Y \hookrightarrow Z \rangle,$$ 
where we again denote by $\overline{X}$ the isomorphism class of $X$. So we have a generating relation of $\K_0\left(\C\right)$ when we have a pushout diagram of the shape
$$\tikz[heighttwo,xscale=2,yscale=2,baseline]{
\node (Y) at (0,1) {$Y$}; 
\node (X) at (0,0) {$X$};
\node (Z) at (1,1) {$Z$};
\node (W) at (1,0) {$X\sqcup_Y Z$};
\node (po)at (0.5,0.5) {$\ulcorner$};
\node at (2,0.5){or};
\node (Y2) at (3,1) {$Y$};
\node (X2) at (3,0) {$X$};
\node (Z2) at (4,1) {$Z$};
\node (W2) at (4,0) {$X\sqcup_Y Z$.};
\node (po)at (3.5,0.5) {$\ulcorner$};
\draw[->]
(Y2) edge (X2)
(Y2)edge (Z2)
(Z2) edge (W2)
(X2)edge  (W2)
(Y) edge (X)
(Y)edge (Z)
(Z) edge (W)
(X)edge  (W);
\draw[right hook ->]
(Y) edge (X)
(Y2)edge (Z2);
}
$$
\end{definition}

\begin{remark}\label{Remark how to dualize pushouts}
When considering the pullback case (i.e. replace pushouts with pullbacks), one needs to replace in the above definition monomorphisms with epimorphisms and in the following proof the cokernel with the kernel.
\end{remark}

We now prove that the above definition indeed is a generalization of the classical definition.

\begin{theorem}\label{Lemma Heapy K0 recovers classical K0}
Let $\A$ be an abelian category. Then we have that $\K_0\left(\A\right)$ recovers the original definition of $\K_0$ by considering the retract along the zero object $\left(\K_0\left(\A\right),+_{\overline{0}}\right)$ i.e. $\overline{A}+_{\overline{0}}\overline{B}=\left[\overline{A},\overline{0},\overline{B}\right]$.
\begin{proof}
Throughout this proof we will distinguish the a priori different constructions by marking the classical constructions with a subscript $c$.

Consider the map 
\begin{align*}
\psi:\K_0\left(\A\right)&\to \K_0\left(\A\right)_c\\
\overline{A}&\mapsto \overline{A}_c.
\end{align*}

We need to first show that this is indeed well defined. Consider for this a generating relation of $\K_0\left(\A\right)$: Let 
\begin{equation}\label{pushout Lemma Heapy K0 revores classical K0}
\tikz[heighttwo,xscale=2,yscale=2,baseline]{
\node (X1) at (0,1) {$B$};
\node (X2) at (0,0) {$A$};
\node (Y1) at (1,1) {$C$};
\node (X+Y) at (1,0) {$D$,};
\node (po)at (0.5,0.5) {$\ulcorner$};

\draw[right hook ->]
(X1) edge node[left] {$\varphi$}  (X2);
\draw[->]
(X1) edge node[above] {$\varphi'$}  (Y1)
(X2) edge  (X+Y)
(Y1)edge (X+Y);
}
\end{equation}
be a pushout inducing $\left[\overline{A},\overline{B},\overline{C}\right]=\overline{D}$. This becomes after passing to $\left(\K_0\left(\A\right),+_{\overline{0}}\right)$ via Lemma~\ref{Proposition passing between groups and heaps with dedicated 0}
\begin{align*}
\left[\overline{A},\overline{B},\overline{C}\right]&=\overline{D}\\
\overline{A}-_{\overline{0}} \overline{B}+_{\overline{0}}\overline{C}&=\overline{D}\\
\overline{A} +_{\overline{0}} \overline{C} &=\overline{B} +_{\overline{0}} \overline{D}
\end{align*} 

 Since $A$ is abelian the pushout (\ref{pushout Lemma Heapy K0 revores classical K0}) corresponds to a short exact sequence 
$$0\to B \xhookrightarrow{\varphi + \varphi'} A\oplus C \twoheadrightarrow D\to 0,$$

which gives rise to the following relation in $\K_0\left(\A\right)_c$
\begin{align*}
\psi\overline{A} + \psi\overline{C}=\overline{A}_c + \overline{C}_c&=\overline{B}_c+\overline{D}_c=\psi\overline{B}+\psi\overline{D}
\end{align*}
and so $\psi$ is well defined.

Furthermore $\psi$ defines a group morphism as we have
\begin{align*}
\psi\left(\overline{A} +_{\overline{0}} \overline{B}\right)&=\psi\left(\left[\overline{A},\overline{0},\overline{B}\right]\right)&*+_{\overline{0}}* = [*,\overline{0},*]\\
&=\psi\left({\overline{A\sqcup_{0}B}}\right) &[*,*,*]=*\sqcup_* *\\
&=\psi\left({\overline{A\oplus B}}\right) &\text{pushout along 0 is } \oplus \\
&=\overline{A\oplus B}_c &\psi\left(\overline{A}\right)=\overline{A}_c\\
&=\overline{A}_c + \overline{B}_c &\text{Lemma~\ref{Lemma split K0 surjects onto ordinary K0}}
\end{align*}

By construction $\psi$ is surjective, in particular we now need to show injectivity. 

To prove $\ker\left(\psi\right)=0$ it is enough to show that for a short exact sequence $0\to A \hookrightarrow B \twoheadrightarrow C \to 0$ we have $\overline{B}=\overline{A} +_{\overline{0}}\overline{C}$. Note that this is equivalent to
\begin{align*}
 \overline{B}&=\overline{A}  +_{\overline{0}} \overline{C}\\
\overline{B}-_{\overline{0}} \overline{A}&=\overline{C}\\
\overline{B}-_{\overline{0}} \overline{A} +_{\overline{0}} 0 &=\overline{C} \iff
\left[\overline{B},\overline{A},\overline{0}\right]=\overline{C}.& \text{Lemma~\ref{Proposition passing between groups and heaps with dedicated 0}}
\end{align*}
But this holds as every short exact sequence is a pushout with a monomorphic leg 
$$\tikz[heighttwo,xscale=2,yscale=2,baseline]{
\node (X1) at (0,1) {$A$};
\node (X2) at (1,1) {$0$};
\node (Y1) at (0,0) {$B$};
\node (X+Y) at (1,0) {$C$,};
\node (po)at (0.5,0.5) {$\ulcorner$};

\draw[->]
(X1) edge  (X2)
(X1)edge  (Y1)
(X2) edge (X+Y)
(Y1)edge (X+Y);
}
$$
and so $\psi$ is injective and hence an isomorphism.
\end{proof}	
\end{theorem}

\begin{remark}
As shown by the proof, Theorem~\ref{Lemma Heapy K0 recovers classical K0} is a direct consequence of the fact that every short exact sequence $0\to A\hookrightarrow B\twoheadrightarrow C \to 0$ in an abelian category $\A$ is a pushout of the form 
$$
\tikz[heighttwo,xscale=2,yscale=2,baseline]{
\node (X1) at (0,1) {$A$};
\node (X2) at (0,0) {$0$};
\node (Y1) at (1,1) {$B$};
\node (X+Y) at (1,0) {$C$,};
\node (po)at (0.5,0.5) {$\ulcorner$};

\draw[->]
(X1) edge  (X2)
(X1)edge  (Y1)
(X2) edge (X+Y)
(Y1)edge (X+Y);
}
$$
and every pushout with a monomorphic leg
$$\tikz[heighttwo,xscale=2,yscale=2,baseline]{
\node (X1) at (0,1) {$B$};
\node (X2) at (0,0) {$A$};
\node (Y1) at (1,1) {$C$};
\node (X+Y) at (1,0) {$D$,};
\node (po)at (0.5,0.5) {$\ulcorner$};

\draw[right hook ->]
(X1) edge node[left] {$\varphi$}  (X2);
\draw[->]
(X1) edge node[above] {$\varphi'$}  (Y1)
(X2) edge  (X+Y)
(Y1)edge (X+Y);
}
$$
in $\A$ is a short exact sequence
$$0 \to A \xhookrightarrow{\varphi + \varphi'} B \oplus C \twoheadrightarrow D \to 0.$$
So the two definitions naturally coincide in this case.
\end{remark}

\begin{proposition}\label{lemma induced morphism of heaps}
Let $F: \C \to \D$ be a functor preserving pushouts and monomorphisms. Then $F$ induces a morphism of heaps:
\begin{align*}
\K_0\left(F\right):\K_0\left(\C\right)&\to \K_0\left(\D\right)\\
\overline{X}&\mapsto \overline{FX}.
\end{align*}
\begin{proof}
As the functor $F$ maps isomorphism classes of objects to isomorphism classes of objects we need to check that $\K_0\left(F\right)$ respects the generating relations. For this consider a pushout giving rise to a generating relation in $\K_0\left(\C\right)$:
$$\tikz[heighttwo,xscale=2,yscale=2,baseline]{
\node (X1) at (0,1) {$Y$};
\node (X2) at (0,0) {$X$};
\node (Y1) at (1,1) {$Z$};
\node (X+Y) at (1,0) {$W$.};
\node (po)at (0.5,0.5) {$\ulcorner$};

\draw[->]
(X1)edge  (Y1)
(X2) edge (X+Y)
(Y1) edge (X+Y);
\draw[right hook ->]
(X1) edge  (X2);
}
$$
Then since $F$ preserves pushouts and monomorphisms we have the pushout 
$$\tikz[heighttwo,xscale=2,yscale=2,baseline]{
\node (X1) at (0,1) {$FY$};
\node (X2) at (0,0) {$FX$};
\node (Y1) at (1,1) {$FZ$};
\node (X+Y) at (1,0) {$FW$.};
\node (po)at (0.5,0.5) {$\ulcorner$};

\draw[->]
(X1)edge  (Y1)
(X2) edge (X+Y)
(Y1) edge (X+Y);
\draw[right hook ->]
(X1) edge  (X2);
}
$$
Which by definition means that $[\overline{FX},\overline{FY},\overline{FZ}]=\overline{FW}$ as claimed.
\end{proof}
\end{proposition}

\begin{remark}
Although the assumptions on $F$ seem very strong, they still are a strict relaxation from $F$ being exact as assumed in the classical constructions.
\end{remark}

\begin{theorem}\label{Lemma canonical projection from split grothendieck group}
Let $\A$ be an additive category admitting pushout. Then we have a canonical epimorphism $$\K_0^\oplus\left(\A\right)\twoheadrightarrow \K_0\left(\A\right),$$ where we endow $\K_0^\oplus\left(A\right)$ with the canonical heap structure given by $$[\overline{A},\overline{B},\overline{C}]_\oplus:= \overline{A}-\overline{B}+\overline{C},$$
see Example~\ref{example Group}  
\begin{proof}
As in the proof of Lemma~\ref{Lemma split K0 surjects onto ordinary K0} we need to prove that the generating relations $
\overline{A}+\overline{B}=\overline{A\oplus B}$ holds in $\K_0\left(\A\right)$, or equivalently $[\overline{A},0,\overline{B}]=\overline{A\oplus B}$. 
For this consider the following diagram exhibiting $X\oplus Y$  as a pushout
$$\tikz[heighttwo,xscale=2,yscale=2,baseline]{
\node (X1) at (0,1) {$0$};
\node (X2) at (0,0) {$A$};
\node (Y1) at (1,1) {$B$};
\node (X+Y) at (1,0) {$A\oplus B$,};
\node (po)at (0.5,0.5) {$\ulcorner$};

\draw[->]
(X1) edge (X2)
(X1)edge (Y1);
\draw[right hook ->]
(X2) edge node[below] {$\iota_A$} (X+Y)
(Y1)edge node [right] {$\iota_B$} (X+Y);
}
$$
where we denote by $\iota$ the canonical inclusions. This shows that the relations in $\K^\oplus_0\left(\A\right)$ are respected by $\K_0\left(\A\right)$ and so we get the desired canonical morphism.
\end{proof}
\end{theorem}

Theorem~\ref{Lemma canonical projection from split grothendieck group} and Theorem~\ref{Lemma Heapy K0 recovers classical K0} show that the heap interpretation of $\K_0\left(\C\right)$ is a natural extension of the classical definition of $\K_0$.

\subsection{Truss structure on $\K_0$ for monoidal categories}

In case the category is monoidal the Grothendieck heap inherits a multiplication making it into a truss, see Definition~\ref{Definition truss}.

\begin{theorem}\label{Theorem Truss structure on K0}
Let $\C$ be a monoidal category admitting pushouts such that $-\otimes -$ preserves pushouts and monomorphisms. Then $\K_0\left(\C\right)$ admits a natural truss structure induced by $$\overline{A}\;\overline{B}=\overline{A\otimes B}.$$
\begin{proof}
As in the case for the classical $\K_0$ we have that $-\otimes -$ induces a natural truss structure on the free heap generated by isomorphism classes $\langle \overline{A}\rangle_{A\in \C}$ given by $\overline{A}\;\overline{B}=\overline{A\otimes B}$. In particular it suffices to prove that the defining relations give rise to an ideal instead of a normal subheap. For this we need to prove that if $\left[\overline{A},\overline{B},\overline{C}\right]=\overline{D}$ is a generating relation, then so are
\begin{align*}
\left[\overline{X}\;\overline{A},\overline{X}\;\overline{B},\overline{X}\;\overline{C}\right]&=\overline{X}\;\overline{D}\\
\left[\overline{A}\;\overline{X},\overline{B}\;\overline{X},\overline{C}\;\overline{X}\right]&=\overline{D}\;\overline{X} 
\end{align*}
for all $X \in \C$.

We will consider the case of left multiplication, the case for right multiplication is analogous.
Let 
$$\tikz[heighttwo,xscale=2,yscale=2,baseline]{
\node (X1) at (0,1) {$B$};
\node (X2) at (0,0) {$A$};
\node (Y1) at (1,1) {$C$};
\node (X+Y) at (1,0) {$D$,};
\node (po)at (0.5,0.5) {$\ulcorner$};

\draw[->]
(X1)edge  (Y1)
(X2) edge (X+Y)
(Y1)edge (X+Y);
\draw[right hook ->]
(X1) edge  (X2);
}
$$

be a pushout inducing $\left[\overline{A},\overline{B},\overline{C}\right]=\overline{D}$. Then, since $ X \otimes -$ preserves pushouts and monomorphisms, we get the pushout 

$$\tikz[heighttwo,xscale=2,yscale=2,baseline]{
\node (X1) at (0,1) {$X\otimes B$};
\node (X2) at (0,0) {$X\otimes A$};
\node (Y1) at (1,1) {$X\otimes C$};
\node (X+Y) at (1,0) {$X\otimes D$,};
\node (po)at (0.5,0.5) {$\ulcorner$};

\draw[->]
(X1)edge  (Y1)
(X2) edge (X+Y)
(Y1)edge (X+Y);
\draw[right hook ->]
(X1) edge  (X2);
}
$$

which gives rise to the relation
\begin{align*}
\left[\overline{X\otimes A},\overline{X\otimes B},\overline{X\otimes C}\right]&=\overline{X\otimes D}\\
\left[\overline{X}\;\overline{A},\overline{X}\;\overline{B},\overline{X}\;\overline{C}\right]&=\overline{X}\;\overline{D}.
\end{align*}
And so the subheap generated by the relations $\left[\overline{A},\overline{B},\overline{C}\right]=\overline{A\sqcup_B C}$ indeed defines an ideal. So we get a canonical truss structure on $\K_0\left(\C\right)$.
\end{proof}
\end{theorem}

\begin{proposition}
Let $\C,\D$ be monoidal categories admitting pushouts such that $-\otimes_\C-$ and $-\otimes_\D-$ preserve pushouts and monomorphisms and let $F:\C \to \D$ be a monoidal functor preserving pushouts and monomorphisms. Then $F$ induces a morphism of trusses.
\begin{align*}
\K_0\left(\C\right) &\to \K_0\left(\D\right)\\
\overline{A} &\mapsto \overline{F\left(A\right)}.
\end{align*}
\begin{proof}
By Theorem~\ref{lemma induced morphism of heaps} we have a canonical morphism of heaps
\begin{align*}
\K_0\left(\C\right) &\to \K_0\left(\D\right)\\
\overline{A} &\mapsto \overline{F\left(A\right)}.
\end{align*}
In particular we need to show that it is compatible with multiplication. For this observe
$$ 
F\left(\overline{A}\;\overline{B}\right)=F\left(\overline{A\otimes_\C B}\right)=\overline{F\left(A\otimes_\C B\right)}=\overline{F\left(A\right)\otimes_\D F\left(B\right)}=\overline{F\left(A\right)}\; \overline{F\left(B\right)}
$$
and so the induced morphism is a morphism of trusses.
\end{proof}
\end{proposition}

\begin{remark}
An essential property for the above constructions was the existence of monomorphisms (respectively epimorphisms) in order to prevent Eilenberg swindle. As monomorphisms and epimorphisms split in triangulated categories they tend to be rather rare and we need to consider a different setup. Luckily a defining property of stable $\infty$-categories is that every cartesian square is cocartesian and vice versa. In particular one can setup similar constructions for the case of stable $\infty$-categories and waldhausen categories.

The machinery and constructions appearing in the suggested construction are significantly deeper results in the theory of homotopy categories. Hence we decided to focus in this paper on the more classical constructions, allowing for a broader audience.

 We are planning to cover this in a future paper. 
\end{remark}

\section{Examples}

\begin{example}\label{Example Sets}

Consider the category of finite sets $\underline{\mathrm{set}}$, this is a category admitting pushouts and a product $\times$ such that $\times$ preserves monomorphisms and pushouts. Then we have the isomorphism of heaps 
$$\K_0\left(\underline{\mathrm{set}}\right)\cong \Z.$$

\begin{proof}
First observe that a finite set $A$ is up to isomorphism completely characterized by its cardinality, so we have a canonical identification $\overline{A}\sim |A|$. Now consider a pushout with a monomorphic leg:
$$\tikz[heighttwo,xscale=2,yscale=2,baseline]{
\node (X1) at (0,1) {$B$};
\node (X2) at (0,0) {$A$};
\node (Y1) at (1,1) {$C$};
\node (X+Y) at (1,0) {$A\sqcup_B C$.};
\node (po)at (0.5,0.5) {$\ulcorner$};

\draw[->]
(X1)edge node[above] {$\varphi$} (Y1)
(X2) edge (X+Y)
(Y1)edge (X+Y);
\draw[right hook ->]
(X1) edge node[right] {$\psi$}  (X2);
}
$$
Then we know that $A\sqcup_B C= A\sqcup C/\sim$  where $a \sim c$ if and only if there is $b\in B$ such that $\psi\left(b\right)=a$ and $\varphi\left(b\right)=c$. Since this relation naturally comes in pairs and $\psi$ is a monomorphism we have that we identify precisely $|B|$ pairs and so $| A \sqcup_B C| = | A | - | B | + | C |$. So 
$$\K_0\left(\underline{\mathrm{set}}\right) \cong \langle \N\mid \left[a,b,c\right]=a -b + c\rangle.$$
This is precisely the heap associated with $\left(\Z,+\right)$.

Now consider the product $- \times -$. Then we have $|A \times B| = |A | | B|$. So we have that the multiplicative structure on $\K_0\left(\underline{\mathrm{set}}\right)$ coincides with the truss structure associated with the ring $\left(\Z,+,\cdot\right)$ and we get that the integers are the decategorification of finite sets.

\end{proof}
\end{example}

\begin{example}\label{Example Top}
Consider the category of topological spaces $\Top$ and let $X$ be homeomorphic to a finite dimensional CW complex . Then we have for some $n\in \N$ the following identity in $\K_0\left(\Top\right)$
\begin{align*}
\overline{X}&=\left[\overline{\sqcup_{i\in I_n} \DDD^n},\overline{\sqcup_{i\in I_n} \Sphere^n},\overline{\sqcup_{i\in I_{n-1}} \DDD^{n-1}},\overline{\sqcup_{i\in I_{n-1}} \Sphere^{n-1}},...,\overline{\sqcup_{i\in I_1} \DDD^1}, \overline{\sqcup_{i\in I_1}\Sphere^1}, \overline{*}\right]\\
\end{align*}
where $\DDD^n$ and $\Sphere^n$ denote the $n$-dimensional disk, respectively sphere.

In particular we get that the Grothendieck heap of topological spaces with CW-type $\K_0\left(\Top_{\text{CW}}\right)$ is generated by the classes of disks and spheres.
\begin{proof}
By assumption $X$ can be built inductively by pushouts of the shape 
$$\tikz[heighttwo,xscale=2,yscale=2,baseline]{
\node (X1) at (0,1) {$\sqcup_{i\in I_k} \Sphere^k$};
\node (X2) at (0,0) {$\sqcup_{i\in I_k} \DDD^k$};
\node (Y1) at (1,1) {$X_{k-1}$};
\node (X+Y) at (1,0) {$X_{k}$};
\node (po)at (0.5,0.5) {$\ulcorner$};

\draw[->]
(X1)edge  (Y1);
\draw[right hook ->]
(X1) edge  (X2)
(X2) edge (X+Y)
(Y1)edge (X+Y);
}
$$
such that there exists an $n\in \N$ with $X_n\cong X$. In particular we have that in $\K_0\left(\Top\right)$ the class of a CW-complex can be identified with the iterated product of its cells. Finally we also get that $\K_0\left(\Top\right)$ admits the structure of a truss as the product of topological spaces respects pushouts and monomorphisms as it is induced by the set theoretic product. 

Altogether we get that the heap $\K_0\left(\Top_{\text{CW}}\right)$ of the category of topological spaces of CW-type is isomorphic to the heap generated by $\Sphere^n$ and $\DDD^n$ together with relations induced by CW-representations.
\end{proof}
\end{example}

 \bibliographystyle{hep}
\bibliography{algebraic_K0_for_unpointed_Categories}  

\end{document}